\documentclass[11pt]{amsart}
\input colordvi

\usepackage{amssymb}
\usepackage{amsmath}

\numberwithin{equation}{section}

\advance\textwidth by 2.5cm 
\advance\oddsidemargin by -1.0cm \advance\evensidemargin by -1.0cm

\newcommand\del[1]{}
\newcommand\think[1]{}
\newcommand\new[1]{}
\newcommand\zus[1]{}

\newcommand\comd[1]{} 
\newcommand\Redd[1]{} 

\def\bdm{\begin{displaymath}}
\def\edm{\end{displaymath}}
\def\bea{\begin{eqnarray}}
\def\eea{\end{eqnarray}}

\newtheorem{theorem}{Theorem}[section]

\newtheorem{defn}[theorem]{Definition}

\newtheorem{Rem}{Remark}

\begin{document}

\date{\today}

\title[Fractional Stochastic active scalar Equations]
{Fractional stochastic active scalar equations generalizing the multi-D-Quasi-Geostrophic \&
2D-Navier-Stokes equations.\\-Short note-}

\author[L.{} Debbi]{Latifa Debbi.\\ {}\\  \SMALL This work is supported by Alexander von Humboldt Foundation.}

\address{Fakult\"at f\"ur Mathematik,
Universit\"at Bielefeld, Universit\"atsstrasse 25,
33615 Bielefeld, Germany.}

\email{ldebbi@math.univ-bielefeld.de}

\maketitle
\vspace{-1.0cm}

\begin{abstract}
We prove the well posedness: global existence, uniqueness and regularity of the solutions, of a class of d-dimensional fractional  stochastic
active scalar equations. This class includes the stochastic, dD-quasi-geostrophic equation, $ d\geq 1$, fractional Burgers equation 
on the circle, fractional nonlocal transport equation and the 2D-fractional vorticity Navier-Stokes equation. We consider the multiplicative noise with locally Lipschitz diffusion term  in
both, the free and no free divergence modes.  The random noise is given by an $Q-$Wiener process with the covariance $Q$ being either of
finite or infinite trace.
In particular, we prove the
existence and uniqueness of a global mild solution for the free divergence
mode in the subcritical regime ($\alpha>\alpha_0(d)\geq 1$), martingale solutions in the general regime ($\alpha\in (0, 2)$) and
free divergence mode, and a local mild solution
for the general mode and subcritical regime. Different kinds of regularity are also established for these solutions. 

\noindent The method used here is also valid for other equations like fractional stochastic velocity Navier-Stokes equations ( work is in progress ). 
The full paper will be published  in Arxiv after a sufficient progress for these equations. 
\vspace{0.15cm}

Keywords: {Scalar active equations, stochastic quasi-geostrophic equation, fractional stochastic vorticity Navier-Stokes equations,
fractional stochastic Burgers equations,
 cylindrical Wiener process, Q-Wiener process, free divergence mode,
martingale solution, mild solution, global and local solutions,
Riesz transform, Hilbert transform, fractional operators, subcritical, critical and supercritical regimes, $\gamma-$radonifying operators,
UMD Banach spaces of type 2.}

\vspace{0.15cm}

Subjclass[2000]: {58J65, 60H15, 35R11.}
\end{abstract}

\section{Introduction}\label{sec-Introduction}

The class of active scalar equations has been introduced as a set of simplistic models for many
complex phenomena in fluid dynamics. Classical examples are the 2D Euler equation in vorticity form and Burgers equation.
Non classical examples are equations with fractional dissipation, such as the quasi-geostrophic equation, the fractional dissipative nonlocal
transport equation and the fractional Burgers equation. But unexpectedly,
it immediately turned out that,
 questions related to this class are not trivial  and far to be achieved.
The fractional scalar equations  showed a great
mathematical and practical potential in the understanding and in the modeling of many complex phenomena, see for short list, \cite{Caffarelli-2009,
Caffarelli-Vasseur2010, SugKak, Pedlosky-87, Sug, Sug-Frac-cal-89}.

\vspace{-0.15cm}


\noindent In this work, we consider a class of d-dimensional active scalar equations driven by a fractional operators and perturbed by
Gaussian random noises.   The representative equation is formally given by,

\begin{equation}\label{Geostro_1Introd}
\left\{
\begin{array}{rl}
 \theta_t (t,x)&= -\nu
(-\Delta)^{\alpha/2}\theta(t,x)+u(t,x)\cdot\nabla \theta(t,x) + G(t, \theta)\xi(t, x), \;\;\\
 u&=(u_j)_j, \;\;\; \theta(0, .)=\theta_0(.),\\
\end{array}
\right.\\
\end{equation}
where, $0< \alpha\leq 2, \; \nu>0$ are real parameters called the fractional order  respectively the viscosity. The spatial coordinate $x\in
\mathbb{T}^d $, with  $ \mathbb{T}^d$  being the d-dimensional Torus, $d\geq 1$ and  $ t\geq 0$ is the time.
The operator $ (-\Delta)^{{\frac\alpha2}}$, denoted bellow by $ A_\alpha$, is the fractional
power of minus Laplacian. The function $\theta_0(.)$ is the initial
condition and  $\xi$ is a random force. The variable $ \theta $ represents the potential temperature and the
 vector field $ u $ represents the fluid velocity. The velocity $ u$ is
determined by the function $\theta $ via the steam function $\psi $. This fact, contrarily to the passive dynamics, makes the dynamics
active. \del{Otherwise the distribution of $\theta $ determines the velocity
 field  $ u $.} In particular, for
\begin{equation}\label{Eq-Exapl-1-u}
 u = (-\frac{\partial \psi}{\partial x_2}, \frac{\partial
\psi}{\partial x_1})= \nabla^\perp \psi,
\end{equation}
Equation \eqref{Geostro_1Introd} is called 2D fractional stochastic vorticity Navier-Stokes equation,
respectively,  stochastic surface quasi-geostrophic, if
\begin{equation}\label{Eq-Exapl-1}
 \theta = \Delta \psi, \;\;\; \text{respectively }\;\;\; (-\Delta)^{\frac12}\psi=
\theta.
\end{equation}
\noindent One possible generalization of the above example, is what called  the modified stochastic quasi geostrophic equation. This latter 
is obtained by considering, see e.g. \cite{Constantin-QG-NS-Arxiv10, Constantin-Iyer-Wu-Global-regularity-Modified-08, 
 Kiselev-active-scalar-11,  Kiselev-Reg-active-scalar-10},
\begin{equation}\label{Eq-Exapl-2}
\theta =  (-\Delta)^{\frac\gamma2}\psi, \;\;\; 1\leq  \gamma \leq 2.
\end{equation}

\noindent In this work, we are dealing with  Equation \eqref{Geostro_1Introd} with nonlinear term being a generalization of
the relation \eqref{Eq-Exapl-2}, for $ d\geq 1, \; \gamma\geq1$ and $ u$ could be or not of divergence free. This equation
covers the fractional stochastic Burgers equation on the circle, the dD stochastic quasi geostrophic equation,
the 2D fractional stochastic vorticity Navier-Stokes equation, d-dimensional stochastic transport equation with nonlocal
coefficients and equations in compressible fluids, such as, the vorticity Navier-Stokes equation with no-zero divergence.

\noindent The novelty in this work, concerns both the equation and the techniques. In fact, to the best knowledge of the author, 
this class of equations is introduced and studied here for the first time, in both deterministic and 
stochastic versions. The idea and the techniques used to prove the global existence of the mild solution for the 
subcritical regime and free divergence mode are also new. 
In particular, we establish a critical threshold,   $1\leq \alpha_0(d)<2$, for which a unique global mild solution exists for this
mode. Other related results are also obtained for larger threshold, $ \alpha_0(d, q)$,  with $ q\geq 2$ is the integrability index. 
The martingale solution is obtained for the 
general regime ($\alpha\in (0, 2)$) and the free divergence mode by generalizing the Hilbert space techniques 
to the Banach space setting.
The local mild solution is proved for the general mode and the subcritical regime by construction. Different kinds of regularity are 
also established for these solutions.

\noindent Let us fix a stochastic basis $ (\Omega, \mathcal{F}, \mathbb{P}, \mathbb{F}, W)$, where
$ (\Omega, \mathcal{F}, \mathbb{P}) $  is a complete probability space, $\mathbb{F} := (\mathcal{F}_t)_{t\geq 0}$
is a filtration satisfying the usual conditions, i.e. $(\mathcal{F}_t)_{t\geq 0}$ is an increasing right continuous filtration. The process
$ W:= (W(t), t\in [0, T])$ is a mean zero
 Gaussian process defined  $ (\Omega,
\mathbb{F}, \mathbb{P}, \mathcal{F} )$, such that the covariance function is given by:
\begin{equation}\label{Eq-Cov-W}
\mathbb{E}[W(t)W(s)]= (t\wedge s)Q,  \;\;\;  \forall\;\;  t,s \geq 0,
\end{equation}
where $ Q $ is a nonnegative operator either of trace class on $ L^2(\mathbb{T}^d)$ or $ Q=I$.
\noindent We rewrite equation \eqref{Geostro_1Introd} as an abstract evolution equation of type,
\begin{equation}\label{Main-stoch-eq}
\Bigg\{
\begin{array}{lr}
 d\theta(t)= \left(-\nu
A_{\alpha}\theta(t) + B^{\sigma, \gamma}(\theta(t))\right)dt+ G(t, \theta(t))dW(t), \; 0< t\leq T,\\
\theta(0)= \theta_0,
\end{array}
\end{equation}
where $ \{\sigma, \gamma\}$ are two parameters characterizing the mode and the regularity of the nonlinear term. In particular, we are interested in
three categories; $ C_a$ for the free divergence mode, $ C_b$ and $ C_c $  for the no free divergence mode with different regularities.

\noindent Let us, before announcing the assumptions on the diffusion term $ G$ and on the initial data $ \theta_0$, 
fix a parameter $\delta$ as follow, 
\begin{itemize}
\item (a) $ \delta\geq 0$, if $ (\sigma, \gamma)\in C_a\cup C_b $.
\item (b) $ \delta > \frac12$, if $ (\sigma, \gamma)\in C_c$.  
\end{itemize}

\noindent {\bf Assumption $(\mathcal{A})$}: We assume that, for $ q\geq 2$  and $ \delta$  given by either  (a) or (b),
the operator $G: H^{\delta, q}(\mathbb{T}^d) \rightarrow \mathcal{L}(L^2, H^{\delta, q})$, is a locally Lipschitz
continuous  and of linear growth map in the following senses, \\
For all $ R>0$, there exists a constant $C_R>0$,  s.t.
\begin{equation}\label{Eq-Cond-Lipschitz-Q-G}
||(G(u)-G(v))Q^\frac{1}{2}||_{R_\gamma(L^{2}, H^{\delta, q})}\leq C_R|u-v|_{H^{\delta, q}},
\;\;\; \forall |u|_{H^{\delta, q}}, |v|_{H^{\delta, q}}\leq R,
\end{equation}
where $ R_\gamma(L^{2}, H^{\delta, q})$ denotes the set of $ \gamma-$radonifying operators from $ L^{2}(\mathbb{T}^d)$ to 
$ H^{\delta, q}(\mathbb{T}^d)$.\\
There exists a constant $ c>0$, s.t
\begin{equation}\label{Eq-Cond-Linear-Q-G}
||G(u)Q^\frac{1}{2}||_{R_\gamma(L^{2}, H^{\delta, q})}\leq c(1+|u|_{H^{\delta, q}}),
\;\;\; \forall u\in H^{\delta, q}(\mathbb{T}^d).
\end{equation}

\noindent {\bf Assumption $(\mathcal{B})$}:  Assume that the initial condition $\theta_0$ is an $\mathcal{F}_0-$random variable satisfying
\begin{equation}\label{Eq-initial-cond}
\theta_0\in L^p(\Omega, \mathcal{F}_0, P; H^{\delta, q_0}(\mathbb{T}^d)),
\end{equation}
with $2\leq q_0< \infty$, $ p\geq 2$ and $ \delta$ is given by either $ (a)$ or $ (b)$ above. 
In the case $ \delta =0$, we allow  $ q_0 =\infty$.

\noindent The full paper is organized as follow, after the introduction, we give in Section 2, the rigorous formulation of Problem 
\eqref{Geostro_1Introd}. 
In particular, we define the nonlinear part of the drift  term,  the different notions of solutions and in the end 
of that section, we announce the main results. In Section 3, we prove several lemmas to estimate 
 the nonlinear term.  Sections  
4-6 are devoted to prove the main theorems of this work; Theorems \ref{Main-theorem-mild-solution-1}-\ref{Main-theorem-mild-solution-3}.


\section{Definitions.}\label{Defs}

\subsection{Definitions of solutions.}\label{Definitions}
In this subsection, we give some definitions of solutions for stochastic partial differential equations. We are interested in solutions as processes in
UMD Banach spaces of type 2. In particular in $ H^{\delta, q}(\mathbb{T}^d), 2\leq q<\infty, \delta\geq 0$. In the cases where such generality is not needed,
as in Definitions \ref{def-variational solution} and \ref{def-martingle-solution} bellow,  we restrict ourselves only on the spaces considered
in this work.  
\begin{defn}\label{def-mild-solution}
Let $ p\geq 2$ and $ X$ be an UMD-Banach space of type 2. Assume that $ \theta_0 \in L^p(\Omega, \mathcal{F}_0, P; X)$. Then  an $ X$-valued adapted
stochastic process $(\theta_t, t\in [0, T])$,  is called mild solution of Equation
\eqref{Main-stoch-eq}, iff
\begin{itemize}
 \item there exists a Banach space, $ X_1\hookleftarrow X$, such that
\begin{equation}\label{eq-cond-cont-mild-solution}
 \theta(\cdot, \omega) \in C([0,T]; X_1)\cap L^\infty(0, T; X),
 \end{equation}
\item  \begin{equation}\label{eq-cond-norm-mild-solution}
 \mathbb{E} \sup_{[0, T]}(|\theta(t)|^p_{X})<\infty
 \end{equation}
\item and $ \forall t \in [0, T]$,
\begin{equation}\label{Eq-Mild-Solution}
\theta(t)= e^{-A_\alpha t}\theta_0 +
\int_0^t e^{-A_\alpha (t-s)}B(\theta (s))ds + \int_0^te^{-A_\alpha (t-s)}G
(\theta(s))W(ds). \; a.s.
\end{equation}
\end{itemize}
\end{defn}

\begin{defn}\label{def-local-mild-solution}
Let $ X$ be an UMD-Banach space of type 2. Assume that $ \theta_0 \in L^p(\Omega, \mathcal{F}_0, P; X)$. A local mild solution of Equation
\eqref{Main-stoch-eq} is a couple $(\theta, \tau_\infty)$, where  $\tau_\infty\leq T\;  a.s. $ is a stopping time and
 $(\theta(t), t\in [0, \tau_\infty))$ is an $ X$-valued adapted
stochastic process such that
\begin{itemize}
\item there exists a Banach space, $ X_1\hookleftarrow X$, such that
\begin{equation}\label{eq-cond-cont-local-mild-solution}
 \theta(\cdot, \omega) \in C([0,\tau_\infty); X_1),
 \end{equation}
 \item  there exists an increasing sequence of stopping time
$ (\tau_n)_{n\in \mathbb{N}}$, s.t. $ \tau_n\nearrow \tau_\infty $.
\item For all $ n\in \mathbb{N}$, the process $ (\theta(t\wedge \tau_n), t\in [0, T])$  satisfies the stopped  \eqref {Eq-Mild-Solution} equation,
i.e. for all $ n \in \mathbb{N}$ and $ \forall t\in [0, T] $,  the following equation is satisfied
\begin{eqnarray}\label{Eq-Mild-Solution-stoped}
\theta(t\wedge \tau_n)= e^{-A_\alpha (t\wedge \tau_n)}\theta_0 &+&
\int_0^{(t\wedge \tau_n)} e^{-A_\alpha (t\wedge \tau_n-s)}B(\theta (s\wedge \tau_n))ds\nonumber \\
& + &\int_0^te^{-A_\alpha (t-s)}1_{[0, \tau_n)}(s)G
(\theta(s\wedge \tau_n))W(ds). \; a.s.
\end{eqnarray}
\end{itemize}
\end{defn}

\begin{defn}\label{def-variational solution}
Let  $ 2\leq q<\infty$. Assume that $ \theta_0 \in L(\Omega, \mathcal{F}_0, P; L^q(\mathbb{T}^d))$.
An $ L^q$-valued adapted stochastic process $(\theta_t, t\in [0, T])$,  is called weak solution of Equation
\eqref{Main-stoch-eq}, iff
 \begin{equation}\label{eq-set-solu-weak}
 \theta(\cdot, \omega) \in L^\infty(0, T; L^q(\mathbb{T}^d))\cap L^2(0, T; H^{\frac \alpha2, 2}(\mathbb{T}^d))\cap
C([0, T]; H^{-\delta', q}(\mathbb{T}^d))\;\;\; a.s.
 \end{equation}
where $ \delta'\geq_1 \max\{\alpha, 1+\frac d{q^*}\}$ 
with  $ q^* = \frac{q}{q-1}$ and such that for all $ \varphi \in D(A_{q^*}^\frac\eta2)$
with $ \eta \geq_2 \max\{1+\frac dq, \frac\alpha2-\frac d2+\frac d{q^*} \}$, we have
\begin{eqnarray}\label{Eq-weak-Solution}
\langle \theta(t), \varphi\rangle &=& \langle \theta_0, \varphi\rangle + \int_0^t \langle \theta (s), A^\frac\alpha2_{q^*}\varphi \rangle ds+
\int_0^t \langle \mathcal{R}^{\gamma, \sigma}\theta (s)\cdot \nabla\varphi, \theta (s)\rangle ds\nonumber\\
& + &\langle\int_0^tG
(\theta(s))dW(s), \varphi\rangle. \; a.s.
\end{eqnarray}
\end{defn}

\begin{defn}\label{def-martingle-solution}
The multiple  $ (\Omega^*, \mathcal{F}^*, \mathbb{P}^*, \mathbb{F}^*, W^*, \theta^*)$, where
$ (\Omega^*, \mathcal{F}^*, \mathbb{P}^*, \mathbb{F}^*, W^*)$ is a stochastic basis with
$ W^* $ being  a $ Q-$Wiener process of trace class and $ \theta^*:=(\theta^*(t), t\in[0, T])$ is an adapted stochastic process,
is called a martingale solution of Equation \eqref{Main-stoch-eq},  iff $ \theta^*$ is a solution of Equation \eqref{Main-stoch-eq} in the sense of 
Definition \ref{def-variational solution} on  the basis $ (\Omega^*, \mathcal{F}^*, \mathbb{P}^*, \mathbb{F}^*, W^*)$.
\end{defn}

\section{Results.}\label{Results}
\noindent The main results of this work are
\begin{theorem}\label{Main-theorem-mild-solution-1}[Free divergence mode \& subcritical regime]
Let  $ d\in\{1,2,3\}$ and $ T>0$ be fixed and let  $ \alpha \in (\alpha_0, 2]$,
with
\begin{equation}\label{Main-value-alpha-0}
\alpha_0 = \alpha_0(d):=  1+ \frac{ d-1}3.
\end{equation}
Assume  $ (\sigma, \gamma) \in C_a$, $ G$  and $ \theta_0 $  satisfying Assumptions  $(\mathcal{A})$
respectively  $(\mathcal{B})$, with $ \max\{2, \frac d{\alpha-1}\} \leq_1 q_0\leq \infty$. Then Equation \eqref{Main-stoch-eq} has a unique global mild solution, $ (\theta(t), t\in [0, T])$,
in the sense of Definition \ref{def-mild-solution}, with $ X= L^{q}(\mathbb{T}^d)$, $ p=q$,
$ X_1= H^{-\delta', q}(\mathbb{T}^d)$, where $ \delta'$ is given in Definition \ref{def-variational solution} and
\begin{itemize}
\item $(1)$  for $ d=1$, then  $ \max\{2, \frac1{\alpha-1}\} \leq_1 q \leq_\infty q_0$,
\item $(2)$  for $ d\in\{2, 3\}$, then $ \frac d{\alpha-1}< q\leq \frac{3d}{d-1}\leq_\infty q_0$.
\end{itemize}
\noindent
The solution also satisfies,
\begin{eqnarray}\label{Eq-reg-theta-mild}
\mathbb{E}\left(\sup_{[0, T]}|\theta(t)|^{q}_{L^{q}}
+\int_0^T|\theta(t)|_{H^{\beta, q}}^2dt \right)<\infty,
\end{eqnarray}
where $ \beta\leq \frac\alpha2-\frac d2+\frac d q$.
If in addition, Assumption $(\mathcal{B}) $ is satisfied for $ \delta >0$, then
\begin{eqnarray}\label{Eq-reg-theta-mild-sobolev}
\theta(\cdot, \omega)\in L^\infty(0, T;\; H^{\delta_1, q}(\mathbb{T}^d)),
\end{eqnarray}
with $ 0\leq \delta_1\leq \min\{ \delta,\; \alpha-1-\frac dq\} $ ($q$ satisfies either $(1)$ or $(2)$).

\noindent Furthermore, for  $ d\in\{2, 3\}$,  the solution exists in the following cases;

{\bf case 1. } $\frac{3d}{d-1} \leq q\leq_\infty \min\{q_0, \frac{2d}{d-\alpha}\}$ and  $ d-2\frac{d}{q} <\alpha \leq 2$.

{\bf case 2. }  $  \frac d{\alpha-1}< q\leq q_0 \leq \frac{3d}{d-1}$ and  $1+\frac{d}{q}< \alpha \leq 2$.

\end{theorem}

\begin{theorem}\label{Main-theorem-mild-solution-2}[General mode \& subcritical regime]
Let $ d \in \mathbb{N}_0$, $ T>0$ and  $ p\geq 2$ be fixed.
Assume $ (\sigma, \gamma) \in C_a\cup C_b\cup C_c$, $ G$  and $ \theta_0 $  
satisfying Assumptions  $(\mathcal{A})$
respectively  $(\mathcal{B})$ with  $ \max\{2, \frac d{\alpha-1}\} \leq_1 q_0\leq \infty$ and that  $ \alpha \in ( 1+ \frac dq,  2]$ with 
$\max\{2, \frac d{\alpha-1}\}\leq_1q\leq_\infty q_0$. Then 
Equation \eqref{Main-stoch-eq} has a local mild solution in the sense of Definition \ref{def-local-mild-solution}, satisfying
\begin{equation}\label{Eq-proper-local-solution}
 \theta  \in L^p(\Omega;  L^\infty(0, \tau_\infty; H^{\delta_1, q}(\mathbb{T}^d))\cap L(\Omega;  C([0, \tau_\infty); H^{-\delta'', q}(\mathbb{T}^d)),
 \end{equation}
with   $\delta_1 $ given in Theorem \ref{Main-theorem-mild-solution-1} and $\delta''\geq \alpha+1+\frac dq-\delta$.
\end{theorem}

\begin{theorem}\label{Main-theorem-mild-solution-3} [General regime \& free divergence mode]
Let $ d \in \mathbb{N}_0$,  $ \alpha \in (0, 2] $ and $ T>0$ be fixed.
Assume $ (\sigma, \gamma) \in C_a$,  $ G$  and   $ \theta_0 $ satisfying  respectively
\eqref{Eq-Cond-Linear-Q-G} and $(\mathcal{B})$. Then  Equation \eqref{Main-stoch-eq}
has a martingale solution, $ (\theta(t), t\in [0, T])$,  in the sense of  Definition \ref{def-martingle-solution},
 satisfying \eqref{eq-set-solu-weak}, \eqref{Eq-weak-Solution} and \eqref{Eq-reg-theta-mild},
with $q=2$.\\ 
\noindent Moreover, the solution $ (\theta(t), t\in [0, T])$ satisfies \eqref{eq-set-solu-weak},  
\eqref{Eq-weak-Solution} and \eqref{Eq-reg-theta-mild}, for  $ d(1-\frac2q)\leq \alpha \leq 2$ and either
\vspace{-0.35cm}
\begin{itemize}
\item $ d\leq \alpha$ and   $ 2\leq q\leq_\infty q_0$, or
\item $ d> \alpha$ and $ 2\leq q \leq_\infty\min\{q_0, \frac{2d}{d-\alpha}\}$.
\end{itemize}
\end{theorem}

\noindent { \bf Acknowledgement.}
The author wishes to thank Professors Ben Goldys, Micheal R\"ockner,  Jan Van Neerven and Dr. Martin Ondrej\`at for
their fruitful discussions. Great thanks also  go to Professor Winfried Sickel for sending
and pointing out important references to the author.





\end{document}